\documentclass[11pt]{article}
\usepackage[a4paper,margin=27mm]{geometry}
\usepackage[T1]{fontenc}
\usepackage{lmodern,microtype}
\usepackage{amsmath,amssymb,amsthm,mathtools}
\usepackage[dvipsnames]{xcolor}
\definecolor{purple}{HTML}{70299B} 
\usepackage[colorlinks=true,linkcolor=MidnightBlue,citecolor=MidnightBlue,urlcolor=MidnightBlue]{hyperref}
\usepackage{enumitem}
\setlist{itemsep=3pt,topsep=5pt}
\newtheorem{theorem}{Theorem}[section]
\newtheorem{lemma}[theorem]{Lemma}
\newtheorem{claim}[theorem]{Claim}
\newtheorem{corollary}[theorem]{Corollary}

\newcommand{\E}{\mathbb E}

\newcommand{\ex}{\operatorname{ex}_k}

\newcommand{\bis}{\operatorname{bis}_w}
\newcommand{\AY}[1]{{#1}}
\newcommand{\GG}[1]{{#1}}

\numberwithin{equation}{section}
\title{Sharp lower bounds for weighted bisections\\of graphs with bounded maximum degree}
\author{
	Stefanie Gerke\thanks{Department of Computing, Security and Mathematics, Royal Holloway University of London, UK,  {\tt stefanie.gerke@rhul.ac.uk}.} \hspace{2mm} Gregory Gutin\thanks{Corresponding author. Department of Computing, Security and Mathematics, Royal Holloway University of London, UK, {\tt g.gutin@rhul.ac.uk}.} \hspace{2mm} Anders Yeo\thanks {Department of Mathematics and Computer Science, University of Southern Denmark, Denmark. {\tt andersyeo@gmail.com}, and Department of Mathematics, University of Johannesburg, South Africa.} \hspace{2mm} Yacong Zhou\thanks{Shenzhen Institutes of Advanced Technology, Chinese Academy of Sciences, China, {\tt yacong.zhou96@gmail.com}.} }

\date{}
\begin{document}
\maketitle

\begin{abstract}
A bisection of a graph is a partition of its vertex set into two parts the sizes of which differ by at most one. We prove that every finite simple
graph of maximum degree at most an odd positive integer $k$, with
arbitrary nonnegative real edge weights, has a bisection of weight at
least $(k+1)/(2k)$ times its total edge weight. The corresponding bound
$(k+2)/(2(k+1))$ for even $k$ follows immediately. Both coefficients
are best possible. Somewhat surprisingly, these are exactly the sharp
coefficients known for maximum weighted cuts, despite the additional
balance requirement. 
\end{abstract}


\section{Introduction}\label{sec:intro}

All graphs considered in this paper are finite and simple, except for
multigraphs obtained explicitly by contraction. Let \(G=(V,E)\) be a
graph with an edge-weight function \(w:E\to\mathbb R_{\ge0}\). Write
$w(G)=\sum_{e\in E}w(e)$. A cut is the set of edges joining the two
parts of a partition of \(V\). A cut is a \emph{bisection} if the part
sizes differ by at most one. Let \(\bis(G)\) denote the maximum weight
of a bisection of \(G\), and let \(\operatorname{cut}_w(G)\) denote the
maximum weight of an unrestricted cut. Thus
\(\bis(G)\le\operatorname{cut}_w(G)\).

The maximum-cut problem and the maximum-bisection problem both seek a
partition with large crossing weight. The balance requirement makes
the latter more restrictive. Pairing the vertices arbitrarily, leaving
at most one vertex unpaired, and \AY{assigning the vertices in each pair to opposite sets independently} and
uniformly gives the elementary bound \(\bis(G)\ge w(G)/2\). A natural
question is how much the maximum-degree bound improves this estimate.

Classical results on cuts in graphs \AY{with} chromatic number \(\chi\)
\cite{AGL,LT,Locke} give a cut with at least a proportion
\[
 p_\chi=\begin{cases}
 \dfrac{\chi+1}{2\chi},&\chi\text{ odd},\\[4pt]
 \dfrac{\chi}{2(\chi-1)},&\chi\text{ even}
 \end{cases}
\]
of all edges, for \(\chi\ge2\). The same argument gives
\(\operatorname{cut}_w(G)\ge p_\chi w(G)\) for arbitrary
nonnegative edge weights. Indeed, fix a proper \(\chi\)-coloring and
choose uniformly a set of \(\lfloor\chi/2\rfloor\) colors for one
part, placing the remaining colors in the other part. Every edge
joins distinct colors and is separated with probability
\[
 \frac{2\lfloor\chi/2\rfloor\lceil\chi/2\rceil}{\chi(\chi-1)}
 =p_\chi.
\]
Linearity of expectation proves the weighted bound. If \(G\) has no
edges, the corresponding weighted statements are immediate.

In particular, a proper \((k+1)\)-coloring of a graph with maximum
degree at most \(k\) yields
\begin{equation}\label{eq:coefficient}
 \operatorname{cut}_w(G)\ge c_k w(G),\qquad
 c_k=\begin{cases}
 \dfrac{k+1}{2k},&k\text{ odd},\\[4pt]
 \dfrac{k+2}{2(k+1)},&k\text{ even}.
 \end{cases}
\end{equation}
These coefficients are sharp for complete graphs with unit weights.
However, balancing the numbers of colors does not balance the numbers
of vertices: color classes can have very different sizes. Thus this
argument \AY{may not give a bisection, so cannot be used to give the same guarantee for bisections.}

The distinction can be substantial even for simple weighted graphs.
For the star \(K_{1,2r-1}\) with unit weights, the maximum cut contains
all \(2r-1\) edges, whereas every bisection has at most \(r\) crossing
edges. In general, therefore, requiring balance can reduce the optimum
to close to half the maximum-cut value. \AY{Therefore,} it is
somewhat surprising that the sharp universal maximum-degree bounds
\eqref{eq:coefficient}, first obtained for unrestricted cuts, continue
to hold for maximum weighted bisections. 

The unweighted case was recently resolved by Guo and Gutin~\cite{GG} improving on a bound by Lee, Loh and Sudakov~\cite{LLS} and a result by  Bollob\'as and Scott~\cite{BS} who had established the exact bound
\((k+1)m/(2k)\) for \(k\)-regular graphs when \(k\) is odd. Their proof cannot be extended to the weighted case. 


Our main result resolves this case.

\begin{theorem}\label{thm:main}
Let \(k\) be an odd positive integer. Every graph \(G\) with
\(\Delta(G)\le k\) and nonnegative real edge weights satisfies
\[
 \bis(G)\ge\frac{k+1}{2k}w(G).
\]
\end{theorem}

Applying Theorem~\ref{thm:main} with \AY{$k+1$ instead of $k$} gives the following
previously known even-degree bound as an immediate consequence.

\begin{corollary}\label{cor:even}
Let \(k\) be a positive even integer. Every graph \(G\) with
\(\Delta(G)\le k\) and nonnegative real edge weights satisfies
\[
 \bis(G)\ge\frac{k+2}{2(k+1)}w(G).
\]
\end{corollary}

Together these statements establish \cite[Conjecture~1.1]{GGYZ} in
full: the balance requirement causes no loss in the sharp universal
weighted-cut coefficients \(c_k\). Both coefficients are best possible: with unit weights,
\(w(K_{k+1})=k(k+1)/2\) and every bisection of \(K_{k+1}\) has weight
\(\lfloor(k+1)^2/4\rfloor\), so
\[
 \frac{\bis(K_{k+1})}{w(K_{k+1})}=c_k.
\]
Neither connectedness nor a positive minimum degree is assumed.

Our proof uses a lower bound based on matchings and in particular Edmonds' characterisation of the matching polytope and its dual. This allows us to establish a probability space and the use of the probabilistic method.

\section{Matching preliminaries}\label{sec:prelim}

For a graph or multigraph \(H\), write \(e(H)=|E(H)|\), counting
parallel edges separately. For disjoint vertex sets
\(X,Y\), let \(e_H(X,Y)\) be the number of edges joining them.
A matching is
\emph{perfect} if it covers every vertex, and \emph{near-perfect} if
it leaves exactly one vertex uncovered. That vertex is called
\emph{exposed}. A nonempty graph \(H\) is \emph{factor-critical} if
\(H-v\) has a perfect matching for every \(v\in V(H)\). A family of
sets is \emph{laminar} if any two of its members are disjoint or one
contains the other.

For the rest of the proof, \(k\) is a fixed odd positive integer.
For every odd vertex set \(T\) in a graph \(H\), \AY{let the excess of $T$ be defined as follows} 
\[
 \ex^H(T)=e(H[T])-\frac{k(|T|-1)}2.
\]

We call \(T\) \emph{overfull} if \(\ex^H(T)>0\), and omit the
superscript when the ambient graph is clear. 
We use Tutte's theorem~\cite{Tutte}:
a graph has a perfect matching if and only if deleting any vertex set
\(A\) leaves at most \(|A|\) odd components.

Note that if $M$ is a matching in $H$ then the incidence vector of $M$ is $(x_e :\ e\in E(H))$ with  $x_e=1$ if $e \in M$ and $x_e=0$, otherwise. 
The first two conditions in Theorem \ref{thm:matchingpolytope} describe all matchings of $H$ provided all $x_e$ are integral. The third condition is needed in order to get the convex hull when we allow
all $x_e$ to be real. 

\begin{theorem}[Edmonds~\cite{Edmonds}]\label{thm:matchingpolytope}
For a loopless multigraph \(H\), the convex hull of the incidence
vectors of its matchings is described by
\begin{equation}\label{eq:matchingpolytope}
 x_e\ge0,\qquad
 \sum_{e\ni v}x_e\le1\quad(v\in V(H)),\qquad
 \sum_{e\in E(H[T])}x_e\le\frac{|T|-1}{2}
 \quad(|T|\text{ odd}).
\end{equation}
\end{theorem}

\begin{corollary}\label{cor:uniformmatching}
Let \(H\) be a loopless multigraph with maximum degree at most \(k\)
and no overfull odd set. For every  real  edge weighting, \(H\)
has a matching of weight at least \(w(H)/k\).
\end{corollary}
\begin{proof}
The vector { $x=(x_e :\ e\in E(H))$ with }$x_e=1/k$ satisfies \eqref{eq:matchingpolytope}.
Its weighted objective is \(w(H)/k\), and some matching has at least
this weight by Theorem~\ref{thm:matchingpolytope}.
\end{proof}

\begin{corollary}\label{cor:nooverfull}
        Let $k$ be an odd integer, let
	\(G\) be a simple graph with maximum degree at most 
	\(k\), and suppose that \(G\) has no overfull odd set. Then, for every
	nonnegative edge weighting,
	\[
	\bis(G)\ge\frac{k+1}{2k}w(G).
	\]
\end{corollary}
\begin{proof}
	By Corollary~\ref{cor:uniformmatching}, there is a matching \(M\) with
	\(w(M)\ge w(G)/k\). The matching-to-bisection bound of
	Gerke, Gutin, Yeo and Zhou~{{\hypersetup{citecolor=purple}\cite[Corollary~2.3]{GGYZ}}} gives
	\[
	\bis(G)\ge\frac{w(G)+w(M)}2
	\ge\frac{k+1}{2k}w(G).
	\]
	For completeness, the first inequality follows by pairing the endpoints
	of each edge of \(M\), pairing the remaining vertices arbitrarily
	(leaving at most one unpaired), and \AY{assigning the vertices in each pair to opposite sets independently} and 
	uniformly. Every matching edge is \AY{in the} cut, and every
	other edge is \AY{in the} cut with probability at least \(1/2\).
\end{proof}

Consequently, any counterexample to Theorem~\ref{thm:main} must contain
an overfull odd set. The remaining arguments control the cost of
expanding such sets after contraction.

\begin{lemma}\label{lem:normalizeddual}
Let \(H\) be a factor-critical graph of odd order
\(h\ge3\), with real edge weights \(w_H\).
Suppose that, for every vertex \(v\), the maximum
weight of a perfect matching of \(H-v\) is the
same number \(c\). Then, there exist
a laminar family \(\mathcal L\) of proper odd vertex sets,
positive numbers \(z_T\) for
\(T\in\mathcal L\), a real number \(\lambda\),
and nonnegative slacks \(s_e\), such that
\begin{equation}\label{eq:dualconstraints}
w_H(e)=\lambda+
\sum_{\substack{T\in\mathcal L\\e\in E(H[T])}}z_T-s_e,
\end{equation}
and
\begin{equation}\label{eq:dualobjective}
c=\frac{h-1}{2}\lambda+
\sum_{T\in\mathcal L}\frac{|T|-1}{2}z_T.
\end{equation}
Moreover, every near-perfect matching of weight
\(c\) contains exactly \((|T|-1)/2\) edges of
\(H[T]\) for each \(T\in\mathcal L\).
\end{lemma}

\begin{proof}
Choose \(C>\sum_{e\in E(H)}|w_H(e)|\) and
replace the edge weight of each edge $e$ by
\[
\widehat{w_H}(e):=w_H(e)+C.
\]
For any two matchings $M,N$, we have
$|w_H(M)-w_H(N)|
\le\sum_{e\in E(H)}|w_H(e)|<C$.
Thus, under the new weights, a matching
with more edges always has larger weight.
In particular, every maximum-weight matching
now has \((h-1)/2\) edges.
By Theorem~\ref{thm:matchingpolytope},
the maximum-weight matching problem is
equivalent to the following linear program:

{\small
\setlength{\jot}{2pt}
\begin{align*}
~~~~~~\text{(P)}\qquad
&\max\quad&
\sum_{e\in E(H)}\widehat{w_H}(e)\cdot x_e\\
& \text{s.t.}\quad
&\sum_{e\ni v}x_e &&\le& \quad 1, 
&& \mbox{for all }v\in V(H),\\
&&\sum_{e\in E_H(T)}x_e&&\le& \quad \frac{|T|-1}{2},
&&  \mbox{for all }T\subseteq V(H),\ |T|\text{ odd},\\
&&x_e&&\ge& \quad 0,
&&  \mbox{for all } e\in E(H).
\end{align*}
}
Thus, its dual is
{\small
\setlength{\jot}{2pt}
\begin{align}
\text{(D)}\qquad
&\min\quad&
\sum_{v\in V(H)}y_v+
\sum_{\substack{T\subseteq V(H)\\|T|\text{ odd}}}
\frac{|T|-1}{2}z_T
\notag\\
&\text{s.t.}\quad&
y_u+y_v+
\sum_{\substack{T\subseteq V(H)\\
|T|\text{ odd}\\uv\in E_H(T)}}z_T
&\ge\widehat{w_H}(uv),
&\quad &\mbox{for all } uv\in E(H),
\label{eq:dualconstraint}\\
&&y_v&\ge0,&\quad &\mbox{for all } v\in V(H),
\notag\\
&&z_T&\ge0,&\quad
&\mbox{for all } T\subseteq V(H),\ |T|\text{ odd}.
\notag
\end{align}
}
Here $y_v$ and $z_T$ are the dual variables
corresponding to the vertex constraints and
the odd-set constraints, respectively.

For every $u\in V(H)$, let $M_u$ be a
maximum-weight perfect matching of $H-u$,
and let $\mathbf{x}_u\in\mathbb{R}^{|E(H)|}$
denote its characteristic vector, i.e.,
$\mathbf{x}_u(e)=1$ if $e\in M_u$ and
$0$ otherwise.
By our assumption, all these matchings have
the same original weight $c$, and hence
their characteristic vectors are optimal
solutions of (P).
\GG{Since every $\mathbf{x}_u$ is optimal for \textup{(P)}, their average
$\mathbf{x}^*=h^{-1}\sum_{u\in V(H)}\mathbf{x}_u$ is optimal for \textup{(P)}.}
Moreover, for every $v\in V(H)$, we have
\[
\sum_{e\ni v}x_e^*
=\frac{1}{h}\sum_{u\in V(H)}
\sum_{e\ni v}\mathbf{x}_u(e)
=\frac{h-1}{h}<1,
\]
where the last equality holds since
$M_u$ contains exactly one edge incident
with $v$ when $u\ne v$, while $M_v$
contains none. Thus, the complementary
slackness theorem implies that every
optimal solution of (D) satisfies
\begin{equation}\label{eq:compslacknessy}
y_v=0 \text{ for all } v\in V(H).
\end{equation}

Among all optimal solutions of (D),
choose one minimizing the potential
\[
\Phi=\sum_{\substack{T\subseteq V(H)\\
|T|\text{ odd}}}
|T|(h-|T|)z_T.
\]
Such a choice is possible since the
optimal dual solutions form a nonempty
polyhedron and $\Phi\ge0$.
Let $\mathcal L'=\{T\subseteq V(H):
|T|\text{ odd},\ z_T>0\}$.

\begin{claim}\label{cl:laminar}
$\mathcal L'$ is a laminar family.
\end{claim}
\begin{proof}
Suppose, for a contradiction, that there
are two crossing sets $A,B\in\mathcal L'$,
i.e., $A\cap B$, $A\setminus B$ and
$B\setminus A$ are all nonempty.
Let $I=A\cap B$ and
$\alpha=\min\{z_A,z_B\}>0$.

First suppose that $|I|$ is even.
Since $A$ and $B$ are odd, both
$A\setminus B$ and $B\setminus A$ are odd.
We decrease $z_A,z_B$ by $\alpha$,
increase $z_{A\setminus B},z_{B\setminus A}$
by $\alpha$, and increase $y_v$ by
$\alpha$ for every $v\in I$.
We claim that every dual constraint
remains satisfied. Indeed, for any edge
$e=uv$, if both endpoints are in $I$,
then the decrease of $z_A$ and $z_B$
is compensated by the increase of
$y_u$ and $y_v$. If both endpoints
are in $A\setminus B$ or both are in
$B\setminus A$, then
the decrease of $z_A$ or $z_B$ is
compensated by the increase of
$z_{A\setminus B}$ or $z_{B\setminus A}$,
respectively.
If one endpoint is in $I$ and
the other is in $A\setminus B$ or
$B\setminus A$, then the decrease of
$z_A$ or $z_B$ is compensated by the
increase of the vertex variable at
the endpoint in $I$.
For all other edges,
the left-hand side of the dual
constraint does not decrease.
Moreover, the dual objective
remains unchanged, since
\[
|I|+\frac{|A\setminus B|-1}{2}
+\frac{|B\setminus A|-1}{2}
=\frac{|A|-1}{2}+\frac{|B|-1}{2}.
\]
Thus, we obtain another optimal
solution of (D) with $y_v>0$ for
every $v\in I$, contradicting
\eqref{eq:compslacknessy}.

Now suppose that $|I|$ is odd.
Then both $I$ and $A\cup B$ are odd.
We decrease $z_A,z_B$ by $\alpha$
and increase $z_I,z_{A\cup B}$
by $\alpha$.
Again, we claim that every dual
constraint remains satisfied.
Indeed, for any edge $e$, if both
endpoints are in $I$, then the
decrease of $z_A,z_B$ is compensated
by the increase of $z_I,z_{A\cup B}$.
If both endpoints are in $A\setminus B$ or
both are in $B\setminus A$, or if one endpoint
is in $I$ and the other is in
$A\setminus B$ or $B\setminus A$, then
the decrease of $z_A$ or $z_B$ is
compensated by the increase of
$z_{A\cup B}$.
If one endpoint is in $A\setminus B$ and
the other is in $B\setminus A$, then the
left-hand side of the dual
constraint increases by $\alpha$
(due to $z_{A\cup B}$).
For all remaining edges, it
remains unchanged.
Moreover, the dual objective remains
unchanged as $|A|+|B|=|I|+|A\cup B|$.
However, the potential $\Phi$ decreases by
\begin{eqnarray*}
&~&\alpha\bigl(
|A|(h-|A|)+|B|(h-|B|)
-|I|(h-|I|)-|A\cup B|(h-|A\cup B|)\bigr)\\
&=&\alpha\bigl(
|A\cup B|^2+|I|^2-|A|^2-|B|^2\bigr)\\
&=&\alpha\bigl(
(|A\setminus B|+|B\setminus A|+|I|)^2
+|I|^2
-(|A\setminus B|+|I|)^2
-(|B\setminus A|+|I|)^2\bigr)\\
&=&2\alpha|A\setminus B||B\setminus A|>0,
\end{eqnarray*}
where the second equality follows from
$|A|+|B|=|I|+|A\cup B|$.
This contradicts the minimality of $\Phi$.
Thus, no two sets in
$\mathcal L'$ cross, and hence
$\mathcal L'$ is laminar.
\end{proof}

Let $\mathcal L=\mathcal L'
\setminus\{V(H)\}$ and
$\lambda=z_{V(H)}-C$.
Note that $z_{V(H)}=0$ if
$V(H)\notin\mathcal L'$.
Thus, for every edge $e\in E(H)$, by
\eqref{eq:dualconstraint},
\eqref{eq:compslacknessy}
and the definition of $\mathcal L'$, we have
\[
\sum_{\substack{T\in\mathcal L'\\
e\in E(H[T])}}z_T\ge w_H(e)+C.
\]
Let $s_e\ge0$ be the slack in this
inequality. Then, by the definitions
of $\mathcal L$, $\lambda$ and
$\widehat{w_H}$, we obtain
\eqref{eq:dualconstraints}.

On the other hand, by strong duality,
the optimal value $p^*$ of (P) equals
the optimal value $q^*$ of (D).
Recall that every maximum-weight
matching under the new weights
$\widehat{w_H}$ is near-perfect,
and the maximum original weight of
a near-perfect matching is $c$.
Thus, $p^*=c+\frac{h-1}{2}C$.
In addition, by \eqref{eq:compslacknessy},
\[
q^*=\sum_{\substack{T\subseteq V(H)\\
|T|\text{ odd}}}
\frac{|T|-1}{2}\cdot z_T
=\sum_{T\in\mathcal L'}
\frac{|T|-1}{2}\cdot z_T,
\]
where the last equality follows
from the definition of $\mathcal L'$.
Thus,
\[
c+\frac{h-1}{2}C=p^*=q^*
=\frac{h-1}{2}z_{V(H)}
+\sum_{T\in\mathcal L}
\frac{|T|-1}{2}z_T.
\]
Substituting $\lambda=z_{V(H)}-C$
gives \eqref{eq:dualobjective}.

Finally, let $M$ be any near-perfect
matching of weight $c$. Then $M$ is
optimal for (P) under the new weights.
Since $z_T>0$ for every $T\in\mathcal L$,
the complementary slackness theorem
implies that the corresponding odd-set
constraint is tight for $M$, i.e.,
\[
|M\cap E(H[T])|=\frac{|T|-1}{2},
\qquad \forall T\in\mathcal L.
\]
This completes the proof.
\end{proof}

\section{An extremal packing of odd sets}\label{sec:packing}

In this section, we always assume that $k$ is
a positive odd integer and $G=(V(G),E(G),w_G)$ is a weighted graph with maximum degree at most $k$.
A family $\mathcal P$ of pairwise disjoint
overfull odd vertex sets in $G$ is called
an \emph{extremal packing} if it maximizes
$\sum_{S\in\mathcal P}\ex(S)$ and,
subject to this, minimizes
$\sum_{S\in\mathcal P}|S|$.
The empty family is allowed.

Note that every overfull odd set $S$ satisfies
\[
\frac{|S|(|S|-1)}{2}
-\frac{k(|S|-1)}{2}
\ge e(G[S])-\frac{k(|S|-1)}{2}
=\ex(S)>0,
\]
and therefore $|S|>k$. Since both $k$ and $|S|$ are odd, we have
\begin{equation}\label{eq:sizeofs}
	|S|\geq k+2.
\end{equation}

\begin{lemma}\label{lem:packing}
	Let $\mathcal P$ be an extremal packing in $G$.
	Then every $G[S]$, $S\in\mathcal P$,
	is factor-critical. Moreover:
	\begin{enumerate}[label=(\roman*)]
		\item for every $S\in\mathcal P$ and every
		family of pairwise disjoint overfull odd
		subsets $T_j\subseteq S$,
		\begin{equation}\label{eq:packingproperty}
			\sum_j\ex(T_j)\le\ex(S);
		\end{equation}
		\item contracting all members of $\mathcal P$,
		and retaining each original external edge
		separately, gives a loopless multigraph $Q$
		with maximum degree at most $k$ and no
		overfull odd set.
	\end{enumerate}
\end{lemma}

\begin{proof}
	Fix an arbitrary $S\in\mathcal P$. We first show that $G[S]$ is factor-critical. Suppose, for a contradiction, that there
	exists a vertex $r\in S$ such that $G[S]-r$
	has no perfect matching.
	Then, by Tutte's theorem,
	there is a subset
	$B\subseteq S\setminus\{r\}$ such that
	$G[S]-(B\cup\{r\})$ has $q>|B|$
	odd components $C_1,\dots,C_q$.
	Since $|S|-1$ is even, $q$ and $|B|$
	\AY{must have the} same parity, and hence
	$q\ge |B|+2$.
	Let $A=B\cup\{r\}$. Then
	\begin{equation}\label{eq:boundonq}
		q\ge |A|+1.
	\end{equation}
	Let $R$ be the union of all even
	components of $G[S]-A$.
	Thus,
    {\small
\setlength{\jot}{2pt}
	\begin{eqnarray*}
		\ex(S)-\sum_{i=1}^q\ex(C_i)
		&=&e(G[S])-\sum_{i=1}^qe(G[C_i])
		-\left(\frac{k(|S|-1)}{2}
		-\sum_{i=1}^q\frac{k(|C_i|-1)}{2}\right)\\
		&=&e(G[S])-\sum_{i=1}^qe(G[C_i])
		-\frac{k(|S|-\sum_{i=1}^q|C_i|+q-1)}{2}\\
		&=&e(G[S])-\sum_{i=1}^qe(G[C_i])
		-\frac{k(|A|+|R|+q-1)}{2}\\
		&\leq&e(G[S])-\sum_{i=1}^qe(G[C_i])
		-k|A|-\frac{k|R|}{2}\\
		&\leq&0,
	\end{eqnarray*}}
	where the second last inequality
	follows from \eqref{eq:boundonq}
	and the last inequality holds since
	$A$ separates the components of $G[S]-A$,
	and therefore all edges between $A$ and
	the $C_i$'s or $R$, as well as edges within
	$A$, are counted in the term $k|A|$,
	while the edges within $R$ are counted
	in the term $k|R|/2$.
	Thus, $\ex(S)\leq\sum_{i=1}^q\ex(C_i)$
	and $\sum_{i=1}^q|C_i|<|S|$
	(as clearly $|A|>0$). Replacing $S$ in $\mathcal P$ by those $C_i$ with $\ex(C_i)>0$ therefore does not
	decrease the total excess and strictly
	decreases the total size,
	a contradiction to the fact that
	$\mathcal P$ is an extremal packing.
	Thus, $G[S]$ is factor-critical.
	
	Property (i) follows directly from the
	choice of $\mathcal P$.
	Indeed, if $\sum_j\ex(T_j)>\ex(S)$,
	we can replace $S$ with the sets $T_j$
	in $\mathcal P$ to obtain a new family
	$\mathcal P'$ with
	\[
	\sum_{U\in\mathcal P'}\ex(U)
	-\sum_{U\in\mathcal P}\ex(U)
	=\sum_j\ex(T_j)-\ex(S)>0,
	\]
	a contradiction to the maximality of
	$\sum_{S\in\mathcal P}\ex(S)$.
	
	It remains to show that (ii) holds.
	Consider a vertex $s\in V(Q)$
	obtained by contracting a set
	$S\in\mathcal P$.
	Then, $d_Q(s)$ is equal to
	\[
	e_G(S,V(G)\setminus S)
	=\sum_{v\in S}
	(d_G(v)-d_{G[S]}(v))
	\leq\sum_{v\in S}
	(k-d_{G[S]}(v))
	=k|S|-2e(G[S])
	=k-2\ex(S),
	\]
	which is less than $k$. In addition, the degrees of all vertices not obtained by contraction remain
	unchanged. Therefore, $\Delta(Q)\le k$.
	
	Now we show that $Q$ has no overfull
	odd sets. Suppose for a contradiction
	that there is an overfull odd set
	$T\subseteq V(Q)$, and
	let $U\subseteq V(G)$ be its
	corresponding vertex set in $G$.
	As all contracted sets have odd order,
	$|U|$ is odd. Moreover, $U$ is an overfull odd set
	in $G$ as
	{\small
\setlength{\jot}{2pt}\begin{eqnarray*}
		\ex^G(U)
		&=&e(G[U])-\frac{k(|U|-1)}{2}\\
		&=&
		e(Q[T])+
		\sum_{\substack{S\in\mathcal P\\S\subseteq U}}
		e(G[S])
		-\frac{k}{2}\left(
		|T|-1+
		\sum_{\substack{S\in\mathcal P\\S\subseteq U}}
		(|S|-1)\right)\\
		&=&\ex^Q(T)+
		\sum_{\substack{S\in\mathcal P\\S\subseteq U}}
		\ex^G(S)\\
		&>&\sum_{\substack{S\in\mathcal P\\S\subseteq U}}
		\ex^G(S),
	\end{eqnarray*}}
	where the last inequality holds as $T$
	is an overfull odd set in $Q$
	(and therefore $\ex^Q(T)>0$).
	
	Thus, replacing all sets
	$S\in\mathcal P$ with $S\subseteq U$
	by $U$ gives us a new family
	$\mathcal P'$ with
	\[
	\sum_{S\in\mathcal P'}\ex(S)
	-\sum_{S\in\mathcal P}\ex(S)
	=\ex(U)-
	\sum_{\substack{S\in\mathcal P\\S\subseteq U}}
	\ex^G(S)>0,
	\]
	a contradiction to the maximality of
	$\sum_{S\in\mathcal P}\ex(S)$.
	This completes the proof of (ii)
	and therefore that of this lemma.
\end{proof}

\section{Local bisections and boundary costs}\label{sec:local}

A near-perfect matching exposing \(v\) is
\emph{$v$-optimal} if its weight is maximum
among all perfect matchings of \(H-v\).

\begin{lemma}\label{lem:local}
	Let $k$ be a positive odd integer,
	$G=(V(G),E(G),w_G)$ a graph \GG{with nonnegative real weights and}
	with maximum degree at most $k$,
	and $\mathcal P$ an extremal packing in $G$.
	For every $S\in\mathcal P$ and $v\in S$,
	let $a_v$ be the maximum weight of a perfect
	matching of $G[S]-v$.
	Then, there exist $b_{G[S]}\ge\max_{v\in S}a_v$
	and a random $\pm1$ signing $(Z_v)_{v\in S}$
	with $\sum_{v\in S}Z_v=1$
	such that
	\begin{equation}\label{eq:localstate}
		2\E[w_G(X_S,Y_S)]
		-w_G(G[S])
		=b_{G[S]}
		\quad \text{and} \quad
		\E[Z_v]\ge0\quad(\text{for all }v\in S),
	\end{equation}
	where $X_S= \{v\in S: Z_v=1\}$ and $Y_S=S\setminus X_S$.
	The nonnegative penalties $t_v=b_{G[S]}-a_v$
	satisfy
	\begin{equation}\label{eq:penalty}
		w_G(G[S])+
		\sum_{v\in S}(k-d_{G[S]}(v))t_v
		\le kb_{G[S]}.
	\end{equation}
\end{lemma}

\begin{proof}
	Fix $S\in\mathcal P$.
	By Lemma~\ref{lem:packing},
	$G[S]$ is factor-critical.
	Let $\nu=\max_{v\in S}a_v$ and
	$q_v=\nu-a_v\ge0$.
	Define a graph $H$ with $V(H)=S$,
	$E(H)=E(G[S])$ and weights
	\[
	w_H(uv):=w_{G[S]}(uv)-q_u-q_v.
	\]
	Thus, for every $v\in S$, the maximum weight
	of a perfect matching of $H-v$ is
	\[
	a_v-\sum_{u\ne v}q_u
	=(a_v-\nu)+\nu-\sum_{u\ne v}q_u
	=\nu-\sum_u q_u,
	\]
	which is the same constant
	$c:=\nu-\sum_uq_u$.
	Thus, applying Lemma~\ref{lem:normalizeddual}
	to $H$, we obtain $\mathcal L$, $z_T>0$,
	$\lambda$, and $s_e\ge0$ satisfying
	\eqref{eq:dualconstraints}--\eqref{eq:dualobjective}.
	Subtracting $k$ times
	\eqref{eq:dualobjective} from the sum of
	\eqref{eq:dualconstraints} gives
	\begin{eqnarray}\label{eq:dualexcess}
		w_H(H)-kc
		&=&\left(e(H)-\frac{k(|S|-1)}{2}\right)\lambda
		+\sum_{T\in\mathcal L}
		\left(e(G[T])-\frac{k(|T|-1)}{2}\right)z_T
		-\sum_{e\in E(H)}s_e
		\notag\\
		&=&\ex(S)\lambda
		+\sum_{T\in\mathcal L}\ex(T)z_T
		-\sum_{e\in E(H)}s_e.
	\end{eqnarray}
	
	We call an undirected edge $rx\in E(H)$
	\emph{free} if there exists a vertex
	$y\in S\setminus\{r,x\}$ such that some
	$r$-optimal matching 
	contains
	$xy$ but $ry\notin E(H)$.
	Let
	\[
	\gamma=\max\bigl(\{0\}\cup
	\{w_H(e):e\text{ is free}\}\bigr).
	\]
	
	\begin{claim}\label{cl:distinctfreeedges}
		For every overfull
		$T\in\mathcal L\cup\{S\}$,
		there are at least $\ex(T)$ distinct
		free edges in $G[T]$.
	\end{claim}
	\begin{proof}
		Note that the number of vertices in $G[T]$
		with degree $k$ is at least
		\[
		\sum_{v\in T}(d_{G[T]}(v)-(k-1))
		=2e(G[T])-k|T|+|T|
		=2\ex(T)-k+|T|.
		\]
		We now show that every such vertex is
		incident with at least one free edge.
		Indeed, for every vertex $r$ with
		$d_{G[T]}(r)=k$, consider a $r$-optimal
		matching.
		If $T\in\mathcal L$, then by
		Lemma~\ref{lem:normalizeddual}, this matching
		contains exactly $(|T|-1)/2$ edges in
		$G[T]$, and hence induces a perfect matching
		of $G[T]-r$.
		If $T=S$, the same conclusion holds since
		the matching exposes $r$ in $H$.
		Since $r$ has exactly $k$ neighbors
		in $T$ and $k$ is odd, some matching edge
		$xy$ has exactly one endpoint, say $x$,
		in the neighborhood of $r$.
		Thus, $rx$ is a free edge in $G[T]$.
		
		Since each free edge can be obtained
		from at most two choices of $r$, and
		$|T|\ge k+2$ by \eqref{eq:sizeofs} as $T$ is overfull,
		the number of distinct free edges in
		$G[T]$ is at least
		\[
		\ex(T)+\frac{|T|-k}{2}
		\geq\ex(T)+1,
		\]
		completing the proof of the claim.
	\end{proof}
	
	Now consider the laminar tree consisting
	of all overfull sets in $\mathcal L$,
	together with $S$ as its root.
	For each set $T$ in this tree, define
	recursively, starting from the leaves,
	\[
	\rho(T)=\max\left\{
	\ex(T),\sum_{C\text{ child of }T}\rho(C)
	\right\}.
	\]
	By the recursive definition, $\rho(S)$
	is the sum of the excesses of some
	pairwise disjoint overfull subsets of $S$.
	Thus, Lemma~\ref{lem:packing}(i) gives
	$\rho(S)\le\ex(S)$ and hence
	$\rho(S)=\ex(S)$.
	
	We next choose $\rho(T)$ distinct
	free edges in $G[T]$ for every set
	$T$ in the tree, starting from the leaves.
	It is possible to choose $\rho(L)$ distinct
	free edges for every leaf $L$ by
	Claim \ref{cl:distinctfreeedges},
	as clearly $\rho(L)=\ex(L)$ in this case.
	Suppose the edges have already been chosen
	for the children of $T$.
	Since these children are disjoint, together
	they contribute
	$\sum_{C\text{ child of }T}\rho(C)$
	distinct free edges in $G[T]$.
	If this number is less than $\rho(T)$,
	then $\rho(T)=\ex(T)$, and
	Claim \ref{cl:distinctfreeedges}
	guarantees that there are enough free edges
	in $G[T]$ to complete the selection.
	Thus, at the root $S$, we obtain a set
	$F$ of exactly $\ex(S)$ free edges such that
	\[
	|F\cap E(G[T])|\ge\ex(T)
	\qquad\text{for every }T\in\mathcal L,
	\]
	since the inequality is automatic
	when $T$ is not overfull, i.e., $\ex(T)\leq 0$.
	
	By \eqref{eq:dualconstraints} and the
	choice of $F$, we have
	\begin{align*}
		\sum_{e\in F}s_e
		&=\ex(S)\lambda+
		\sum_{T\in\mathcal L}
		|F\cap E(G[T])|z_T
		-\sum_{e\in F}w_H(e)\\
		&\ge\ex(S)\lambda+
		\sum_{T\in\mathcal L}\ex(T)z_T
		-\ex(S)\gamma.
	\end{align*}
	Since $s_e\ge0$ for all $e\in E(H)$,
	combining this with \eqref{eq:dualexcess}
	gives
	$w_H(H)-kc\le\ex(S)\gamma$.
	Recall that
	$c=\nu-\sum_{v\in S}q_v$ and
	$w_H(H)=w_G(G[S])-
	\sum_{v\in S}d_{G[S]}(v)q_v$.
	Therefore,
	\begin{equation}\label{eq:freebound}
		w_G(G[S])-k\nu+
		\sum_{v\in S}(k-d_{G[S]}(v))q_v
		\le\ex(S)\gamma.
	\end{equation}
	
	Let $b_{G[S]}=\nu+\gamma/2$.
	We now construct a random signing satisfying
	\eqref{eq:localstate}.
	
	If $\gamma=0$, choose a vertex $r\in S$ with $a_r=\nu$
	and a $r$-optimal matching,
	put $r$ on the $+1$ side, and independently assign to each
	matching pair $(Z_x,Z_y)=(+1,-1)$ or
	$(-1,+1)$, each with probability $1/2$.
	Then $\sum_{v\in S}Z_v=1$, the expected
	surplus is $\nu=b_{G[S]}$, and
	$\E[Z_r]=1$ while $\E[Z_v]=0$
	for every $v\ne r$.
	Hence \eqref{eq:localstate} holds.
	
	Otherwise, let $rx$ be a free edge with
	$w_H(rx)=\gamma>0$, witnessed by a
	$r$-optimal matching
	that contains $xy$, where
	$ry\notin E(H)$.
	Consider the following two random signings:
	
	(a) Put $r,y$ on the $+1$ side
	and $x$ on the $-1$ side, and independently assign to each
	remaining matching pair
	$(Z_u,Z_v)=(+1,-1)$ or $(-1,+1)$,
	each with probability $1/2$.
	
	(b) Choose an $x$-optimal matching,
	put $x$ on the $+1$
	side, and independently assign to each
	matching pair $(Z_u,Z_v)=(+1,-1)$ or
	$(-1,+1)$, each with probability $1/2$.
	
	Both signings have sum $1$.
	In (a), all matching edges and $rx$
	are included in the cut, while no edge
	is forced to be excluded since
	$ry\notin E(H)$.
	Thus, its expected surplus is
	$a_r+w_G(rx)$, while in (b) it is $a_x$.
	
	Now choose (a) or (b), each with
	probability $1/2$.
	Since $Z_r=Z_y=1$ and $Z_x=-1$ in (a),
	while $Z_x=1$ and $\E[Z_v]=0$
	for every $v\ne x$ in (b),
	we have $\E[Z_r]=\E[Z_y]=1/2$
	and $\E[Z_v]=0$ for every
	$v\notin\{r,y\}$.
	Thus, $\E[Z_v]\ge0$ for all $v\in S$.
	Moreover, the expected surplus is
	\[
	\frac{a_r+w_G(rx)+a_x}{2}
	=\nu+\frac{w_G(rx)-q_r-q_x}{2}
	=\nu+\frac{\gamma}{2}
	=b_{G[S]}.
	\]
	This completes the proof of \eqref{eq:localstate}.
	
	Finally, since
	$t_v=b_{G[S]}-a_v=q_v+\gamma/2\ge0$,
	it follows from \eqref{eq:freebound} that
	\begin{align*}
		w_G(G[S])+
		\sum_{v\in S}(k-d_{G[S]}(v))t_v
		&\le k\nu+\ex(S)\gamma
		+\frac{\gamma}{2}
		\sum_{v\in S}(k-d_{G[S]}(v))\\
		&=k\nu+\ex(S)\gamma
		+\frac{\gamma}{2}(k-2\ex(S))\\
		&=kb_{G[S]}.
	\end{align*}
	This proves \eqref{eq:penalty}
	and completes the proof.
\end{proof}

\section{Proof of the main theorem}\label{sec:assembly}

Now, we are ready to give the proof of our main result.

\begin{proof}[Proof of Theorem~\ref{thm:main}]
	Let $\mathcal P$ be a family of pairwise disjoint overfull odd vertex sets that maximizes $\sum_{S\in\mathcal P}\ex(S)$ and, subject to this, minimizes $\sum_{S\in\mathcal P}|S|$. For every $S\in\mathcal P$, let $b_{G[S]}$, $a_v$ and $t_v$ ($v\in S$) be as in Lemma~\ref{lem:local}. For every $v\notin\bigcup_{S\in\mathcal P}S$, let $t_v=0$.
	
	Contract all sets in $\mathcal P$ to obtain a multigraph $Q$, retaining all external edges separately. Let $c$ be the natural bijection from the edges of $G$ surviving the contraction to $E(Q)$. For every such edge $e$ with endpoints $u$ and $v$, assign
	\[
	w_Q(c(e)):=w_G(e)-t_u-t_v.
	\]
	Then, by Lemma~\ref{lem:packing}(ii), $Q$ has maximum degree at most $k$ and no overfull odd vertex sets. Thus, by Corollary~\ref{cor:uniformmatching}, $Q$ has a matching $M'$ with
	\begin{equation}\label{eq:quotientmatching}
		w_Q(M')\ge\frac{w_Q(Q)}{k}.
	\end{equation}
	
	Now we lift $M'$ to a matching $M$ of $G$. Let $C\subseteq V(Q)$ be the set of new vertices obtained from the contraction. For every $e\in M'$, let $u$ and $v$ be the endpoints of $c^{-1}(e)$ in $G$, and include $c^{-1}(e)$ in $M$. If $u\in S$ for some $S\in\mathcal P$, then we also include a maximum-weight perfect matching of $G[S]-u$, whose weight is $a_u$. We do the same for $v$ if it belongs to some set in $\mathcal P$. Since $M'$ is a matching, the resulting set $M$ is also a matching of $G$. For any edge $e\in M$ with endpoints $x$ and $y$, 
we independently assign \AY{$(X_x,X_y)$ to be $(+1,-1)$ or $(-1,+1)$ with probability $1/2$.}
	
	It remains to assign signs to the vertices not covered by $M$. Let $\mathcal U$ be the family of vertex sets of $G$ corresponding to the vertices of $Q$ not covered by $M'$. Each $A\in\mathcal U$ is either a set in $\mathcal P$ or a singleton vertex outside all sets in $\mathcal P$. For every $A\in\mathcal U$, we independently choose a random signing $(Z_x)_{x\in A}$ such that $\sum_{x\in A}Z_x=1$ and \AY{$\E(Z_x)\ge 0$} for all $x\in A$: If $A\in\mathcal P$, we use the probability space given by Lemma~\ref{lem:local}; if $A=\{u\}$ is a singleton, we simply set $Z_u=1$.
	
	We pair the sets in $\mathcal U$ arbitrarily so that at most one set is left over. For every paired two sets $A,B\in\mathcal U$, we put their larger partite sets on opposite sides. More precisely, we independently choose a uniform random sign $\sigma\in\{-1,+1\}$ and set $X_x=\sigma Z_x$ for $x\in A$ and $X_y=-\sigma Z_y$ for $y\in B$. If one set $A$ remains unpaired, we choose another independent uniform random sign $\sigma$ and set $X_x=\sigma Z_x$ for all $x\in A$. Since each signing $(Z_x)_{x\in A}$ has sum $1$, the resulting signing $X$ gives a bisection $(X^{-1}(+1),X^{-1}(-1))$ of $G$.
	
	Now consider any external edge $e$ with endpoints $x\in A$ and $y\in B$, where $A,B\in\mathcal U$ are paired together. We claim that such edges are included in this random bisection with probability at least $1/2$. Indeed, it is easy to check that $\frac{1-X_xX_y}{2}$ is the indicator for whether this edge is included in the bisection. Thus,
	\[
	\Pr[e\text{ is included}]
	=\E\left[\frac{1-X_xX_y}{2}\right]
	=\E\left[\frac{1+Z_xZ_y}{2}\right]
	=\frac{1+\E[Z_x]\E[Z_y]}{2}
	\ge\frac12,
	\]
	where the last inequality follows from Lemma~\ref{lem:local} and the definition of $Z$ for singleton sets.
	
	Every edge in $M$ is always included in the bisection. Now consider an external edge $e$ not in $M$. If at least one endpoint of $e$ is covered by $M$, then since we independently assigned to $+1$ or $-1$ to the random variable corresponding to its endpoint in $M$ with probability $1/2$, \AY{it is} easy to show that $e$ is included with probability exactly $1/2$.
	Otherwise, both endpoints of $e$ belong to sets in $\mathcal U$, say $A$ and $B$. If $A$ and $B$ are paired together, we have already shown that $e$ is included with probability at least $1/2$. If not, the signs of its endpoints are independent and uniformly distributed over $\{-1,+1\}$, since the random signs $\sigma$ assigned to different pairs are independent and uniform. Hence, $e$ is included with probability exactly $1/2$. Consequently, every external edge not in $M$ is included with probability at least $1/2$.
	
	Now consider the surplus contributed by the internal edges of the sets in $\mathcal P$ and the external edges in $M$. Let $e$ be an edge of $G$ with $c(e)\in M'$ and endpoints $u\in S_u$ and $v\in S_v$, where $S_u,S_v\in\mathcal P$. Since $e$ and the two chosen internal matchings are always included in the bisection, while all other internal edges of these two sets are included with probability $1/2$, their total contribution to the surplus is
	\[
	w_G(e)+a_u+a_v
	=b_{G[S_u]}+b_{G[S_v]}-t_u-t_v+w_G(e)
	=w_Q(c(e))+b_{G[S_u]}+b_{G[S_v]}.
	\]
	If one or both endpoints of $e$ do not belong to any set in $\mathcal P$, the same calculation holds by omitting the corresponding terms. On the other hand, every set $S\in\mathcal P$ whose contracted vertex is not covered by $M'$ contributes an expected surplus of $b_{G[S]}$ as we used the probability space in Lemma~\ref{lem:local} for it. Since all other external edges contribute a nonnegative expected surplus, summing over all edges in $M'$ and all unmatched sets in $\mathcal P$ gives
	\begin{equation}\label{eq:globalcut}
		2\E [w_G(X^{-1}(+1),X^{-1}(-1))]-w_G(G)
		\ge \sum_{S\in\mathcal P}b_{G[S]}+w_Q(M')
		\ge\sum_{S\in\mathcal P}b_{G[S]}+\frac{w_Q(Q)}{k},
	\end{equation}
	where the last inequality follows from \eqref{eq:quotientmatching}.

	On the other hand, as $w_Q(c(e))=w_G(e)-t_u-t_v$ for every edge $e$ surviving the contraction, where $u$ and $v$ are its endpoints in $G$, we have
	\begin{align}
		w_G(G)
		&=w_Q(Q)+\sum_{S\in\mathcal P}
		\left(w_G(G[S])+
		\sum_{u\in S}e(u,V(G)\setminus S)t_u\right)
		\notag\\
		&\le w_Q(Q)+\sum_{S\in\mathcal P}
		\left(w_G(G[S])+\sum_{u\in S}(k-d_{G[S]}(u))t_u\right)
		\notag\\
		&\le w_Q(Q)+k\sum_{S\in\mathcal P}b_{G[S]},
		\label{eq:globalweight}
	\end{align}
	where the first inequality follows from the fact that $\Delta(G)\leq k$ (and therefore
	$e(\{u\},V(G)\setminus S)\le k-d_{G[S]}$) and $t_u\ge0$, and the second follows from Lemma~\ref{lem:local}.

	Combining \eqref{eq:globalcut} and \eqref{eq:globalweight}, we obtain
	\[
	\E [w_G(X^{-1}(+1),X^{-1}(-1))]
	\ge \frac{w_G(G)}2
	+\frac12\left(\sum_{S\in\mathcal P}b_{G[S]}
	+\frac{w_Q(Q)}k\right)
	\ge \frac{k+1}{2k}w_G(G).
	\]
	Hence there exists a bisection of $G$ with weight at least $(k+1)w_G(G)/(2k)$, as desired.
\end{proof}

\section{Concluding remarks}

Theorem~\ref{thm:main} and Corollary~\ref{cor:even} show that the sharp
maximum-degree coefficients for unrestricted cuts also hold for
bisections with arbitrary nonnegative edge weights. The proof does
not require the graph to be regular or connected.

Corollary~\ref{cor:nooverfull} settles the case without overfull odd
sets directly by combining the matching-polytope theorem with the
matching-to-bisection bound. The rest of the proof reduces the general
case to a quotient with no overfull odd set while controlling the
losses in expanding the contracted blocks.

The local inequality \eqref{eq:penalty} accounts for the loss incurred
when a prescribed vertex must be exposed in a selected block. 
\GG{After assigning each external edge \(uv\) the adjusted weight
\(w_G(uv)-t_u-t_v\), the matching-polytope bound gives a quotient
matching of weight at least \(w_Q(Q)/k\). Lifting this matching and
applying the local boundary-cost inequality completes the proof.}

%

\section*{Declaration on the Use of Generative AI} 
\GG{The authors used ChatGPT 6 Astra to assist in discussing proof
strategies and checking proofs. The authors take full responsibility for the content of the
manuscript.}


\end{document}